\documentclass[12pt]{article}

\usepackage{amsfonts}
\usepackage{xr}
\usepackage{graphicx}
\usepackage{amsmath}

\def\to{\rightarrow}
\def\Z{\mathbb{Z}}
\def\R{\mathbb{R}}
\def\N{\mathbb{N}}

\def\o{\omega}
\def\O{\Omega}

\def\E{\mathbb{E}}
\def\P{\mathbb{P}}

\def\tilde{\widetilde}
\def\epsilon{\varepsilon}
\def\trait (#1) (#2) (#3){\vrule width #1pt height #2pt depth #3pt}
\def\fin{\hfill\trait (0.1) (5) (0) \trait (5) (0.1) (0) \kern-5pt \trait (5) (5) (-4.9) \trait (0.1) (5) (0)}
\def\e{\epsilon}

\newtheorem{thm}{\bf Theorem}[section]
\newtheorem{lem}[thm]{\bf Lemma}
\newtheorem{prop}[thm]{\bf Proposition}
\newtheorem{cor}[thm]{\bf Corollary}

\title{How does the spreading speed associated with the Fisher-KPP equation depend on random stationary diffusion and reaction terms?}
\date{}
\author{Gregoire Nadin\thanks{CNRS, UMR 7598, Laboratoire Jacques-Louis Lions, F-75005 Paris, France}}

\begin{document}
\maketitle

\begin{flushright}
{\em In honor of Chris Cosner's sixtieth birthday.}
\end{flushright}

\begin{abstract}
We consider one-dimensional reaction-diffusion equations of Fisher-KPP type with random stationary ergodic coefficients. 
A classical result of Freidlin and Gartner \cite{Gartner} yields that the solutions of the initial value problems associated 
with compactly supported initial data admit a linear spreading speed almost surely. We use in this paper a new characterization of this spreading speed recently proved in 
\cite{BN1} in order to investigate the dependence of this speed with respect to the heterogeneity of the diffusion and reaction terms. We prove in particular that adding a reaction term with null 
average or rescaling the coefficients by the change of variables $x\to x/L$, with $L>1$, speeds up the propagation. From a modelling point of view, 
these results mean that adding some heterogeneity in the medium gives a higher invasion speed, while fragmentation of the medium slows down the invasion.

\vspace{0.2in}\noindent \textbf{Key words}: Eigenvalue optimization; reaction-diffusion equations; spreading properties; random stationary ergodic coefficients; 
biological invasions. 
\noindent \textbf{Aknowledgements}: The research leading to these results has received funding from the European Research Council 
under the European Union's Seventh Framework Programme (FP/2007-2013) / ERC Grant
Agreement n.321186 - ReaDi -Reaction-Diffusion Equations, Propagation and Modelling.
\end{abstract}


\section{Introduction}

\subsection{Framework and hypotheses}\label{sec:hyp}

We investigate in this paper reaction-diffusion equations with coefficients depending on a space variable $x\in\R$ and on a random event $\o\in\O$:

\begin{equation} \label{Cauchy}
\left\{ \begin{array}{lcl}
\partial_t u - \partial_x \big( a(x,\omega) \partial_x u\big) = f(x,\omega,u) \quad &\hbox{in}& \quad (0,\infty) \times \R\times \O,\\
u(0,x,\omega) = u_0 (x) \quad &\hbox{over}& \quad \R\times \O,\\
\end{array} \right. \end{equation}
where $u_0\not\equiv 0$ is a compactly supported, continuous and nonnegative initial datum. The dependence in $(x,\o)$ means that a space 
heterogeneous equation is associated with any $\o\in\O$. Our aim is to determine the link between the asymptotic behaviour of the associated solution
$(t,x) \mapsto u(t,x,\o)$ and the heterogeneity of the coefficients $a$ and $c$.  

If the coefficients do not depend on $(x,\o)$, we recover the classical Fisher-KPP equation \cite{Fisher, KPP}. This equation is involved in many models 
of genetics \cite{Aronson} or population dynamics \cite{Shigesadabook}. It is natural when investigating biological invasion processes to consider heterogeneous media, that is, heterogeneous
diffusion and reaction terms in the Fisher-KPP model. Periodic Fisher-KPP equations have been extensively studied in the past decade 
\cite{BHNa, BHR2, ElSmaily, ElSmailyKirsch, Freidlin, Gartner, HamelNadinRoques, dependence, RyzhikZlatos, Shigesada1, Shigesadabook, Zlatosc/A}. 
People are now focusing on more general classes of heterogeneities. We consider here random stationary ergodic coefficients. 

To be more specific, we consider a probability space $(\O,\P,\mathcal{F})$ and we assume that the derivative of the reaction rate $(x,\omega)\in \R\times \O\mapsto f_s(x,\o,0)$
and the diffusion coefficient $a:\R\times \O \to (0,\infty)$ are random stationary ergodic variables. That is, 
there exists a group $(\pi_x)_{x\in\R}$ of measure-preserving transformations of $\O$ such that
\begin{equation} a(x+y,\omega) = a(x,\pi_y \omega) \quad \hbox{and} \quad f_s(x+y,\omega,0)=f_s(x,\pi_y \omega,0)
\quad \hbox{for all} \quad (x,y,\omega)\in \R\times\R\times \O\end{equation} 
and if $\pi_x A=A$ for all $x\in\R$ and for a given  $A\in\mathcal{F}$, then $\P(A)=0$ or $1$. 
This hypothesis heuristically means that the statistical properties of the medium does not depend on the place where we observe it. This is a very natural hypothesis. Indeed, periodic and almost periodic 
coefficients could be considered as particular random stationary ergodic coefficients \cite{PapanicolaouVaradhan}. 

We consider reaction terms of KPP type, meaning that for all $x\in\R$ and almost every $\o\in\O$:
\begin{equation} \label{hyp:eq} f(x,\o,0) = f(x,\o,1)= 0, \quad 
\inf_{x\in\R} f(x,\o,s)>0 \hbox{ for all } s\in (0,1),\end{equation}
\begin{equation} \label{hyp:KPP} f(x,\omega,s) \leq c(x,\o)s:= f_s(x,\omega,0)s \hbox{ for all } s\in [0,1].\end{equation} 
A typical nonlinearity satisfying this set of hypotheses is $f(x,\o,s) = c(x,\o) s(1-s)$, where $c$ is a random stationary ergodic variable such that 
$\inf_{x\in\R} c(x,\o)>0$ for almost every $\o\in\O$. 

Lastly, $\inf_{x\in\R} a(x,\o)>0$ for almost every $\o\in\O$ and we require the functions $a$, $a'$, $f$ and $f_s(\cdot,0)$ to be almost surely uniformly
continuous and bounded with respect to $x$ uniformly in $s$ and $\mathcal{C}^{1+\gamma}$ with respect 
to $s$ uniformly in $x$. Here and throughout the paper, $a'$ denotes the derivative of $a=a(x,\o)$ with respect to $x$.


\subsection{Definition of the spreading speed}

We are interested in the spreading properties related to equation (\ref{Cauchy}). That is, we want to determine the location of the level sets of $u(t,\cdot,\o)$
at large times. The pioneering works in this direction have been carried out by Freidlin and Gartner \cite{Gartner} and Freidlin \cite{Freidlin2, Freidlin}, 
through a fully probabilistic approach. These earlier works assumed that the whole linearity $(x,\o)\mapsto f(x,\o,s)$ was random stationary ergodic for all $s\in [0,1]$. 
Spreading properties in the framework described above, where only the linearization $(x,\o)\mapsto f_s(x,\o,0)$ is random stationary ergodic, were 
derived by Berestycki and the author in \cite{BN1}. Let now state these spreading properties rigorously. 

\begin{thm} \cite{Gartner, Freidlin, Freidlin2, BN1}\label{thm:Freidlin}
Under the previous hypotheses, there exists a speed $w^*>0$ such that the solution $u=u(t,x,\o)$ of (\ref{Cauchy}) satisfies for almost every $\o\in\O$:
\begin{equation}\label{eq:spreading} \left\{\begin{array}{lllcl}
\hbox{for all  }w\in (0,w^*),& \ \lim_{t\to+\infty}
\sup_{x\in [0,wt)}& |u(t,x,\o)-1|  &= & 0,\\
\hbox{for all } w> w^*,& \
\lim_{t\to+\infty} \sup_{x\geq wt} &|u(t,x,\o)|  &= & 0.\\
\end{array}\right. \end{equation}
Moreover, the speed $w^*=w^*(a,f)$ is characterized by (\ref{characw}).
\end{thm}

Heuristically, this result means that an observer who moves with speed $w$ will only see the stable steady state $1$ if he moves too slowly ($w<w^*$) and the unstable one 
$0$ if he moves too fast ($w>w^*$). This theorem has been extended to multidimensional space-time random stationary ergodic diffusions, with homogeneous reaction terms, 
by Nolen and Xin \cite{NolenXinrandomgeneral}. 
We also refer here to the works of Souganidis \cite{Souganidis} and Lions and Souganidis \cite{LionsSouganidis} who proved some related homogenization 
results. Let us also mention that spreading properties have recently been extended to one-dimensional reaction-diffusion equations of ignition type 
by Nolen and Ryzhik \cite{NolenRyzhik}. In this case the spreading speed is determined by the global mean speed of random travelling waves (which existence remains 
unclear in the Fisher-KPP framework). 

Let now interpret this result in terms of biological modelling. In reaction-diffusion models, $u$ is the density of a population, $f(x,\o,u)$ is the growth rate, 
$c(x,\o)=f_s(x,\o,0)$ is the growth rate per capita at small density and $a(x,\o)$ is the diffusion rate. Theorem \ref{thm:Freidlin} means that the population will
invade the environment at speed $w^*=w^*(a,f)$. A natural question is thus to determine whether the heterogeneity of the coefficients $a$ and $f$ speeds up or slows down the
invasion. Of course it is difficult to get such dependence results from definition (\ref{eq:spreading}) of the spreading speed $w^*$ and alternative characterizations
of $w^*$ would be helpful. Freidlin \cite{Freidlin} gave such a characterization. 


\subsection{Freidlin's characterization of the spreading speed}

Define
\begin{equation} \label{def:Lambda1}
\Lambda_1^\o= \Lambda_1^\o (a,c) := \sup_{\alpha \in H^1 (\R)\backslash \{0\}} \displaystyle\frac{\int_\R \Big( -a(x,\o)\alpha'(x)^2+ c(x,\o) \alpha (x)^2\Big)dx}{\int_\R \alpha^2 (x)dx}.   
\end{equation}
This quantity is almost surely deterministic, that is, it is almost everywhere identical to a quantity which does not depend on $\o$ (see \cite{Nolen}). We will 
thus forget its dependence with respect to $\o$ and denote it $\Lambda_1 (a,c)$ when there is no ambiguity.  
Then for all $\gamma > \Lambda_1 (a,c)$ and for almost every $\o\in\O$, there exists a unique positive solution $\phi_\gamma (\cdot,\o) \in \mathcal{C}^2 (\R)$ of 
\begin{equation} \left\{\begin{array}{l}
\big( a(x,\o)\phi_\gamma'\big)' +c(x,\o)\phi_\gamma = \gamma\phi_\gamma \hbox{ in } \R\times \O,\\
\phi_\gamma (0,\o)=1, \quad \lim_{x\to +\infty} \phi_\gamma (x,\o) = 0.\\
                 \end{array}\right. \end{equation}
The existence and uniqueness of $\phi_\gamma$ had been proved through probabilistic techniques by Freidlin \cite{Freidlin} and through PDE tools by Nolen \cite{Nolen}. 
This uniqueness yields that $\phi_\gamma(x+y,\o)=\phi_\gamma (x,\pi_y\o)\phi_\gamma (y,\o)$ and thus the Birkhoff ergodic theorem ensures that the limit 
$$\mu (\gamma):= -\lim_{x\to +\infty} \frac{1}{x} \ln \phi_\gamma (x,\o) \quad \hbox{is well-defined and does not depend on } \o\in\O \hbox{ almost surely}.$$
Moreover, this limit is positive and convex with respect to $\gamma>\Lambda_1$.
The characterization of $w^*$ reads \cite{Freidlin}:
\begin{equation}\label{characw}
 w^* = \min_{\gamma>\Lambda_1} \frac{\gamma}{\mu (\gamma)}.
\end{equation}

Such a formula could be useful to investigate the link between the heterogeneity of the coefficients $a$ and $f$ and the spreading speed $w^*$. 
Surprisingly, there are very few results on this topic in the literature. 

The case of multi-dimensional reaction-diffusion equations with a space-time random advection term has been fully investigated by Nolen and Xin at the end of the 2000's. 
They proved that there exists a unique spreading speed in this space-time heterogeneous framework \cite{NolenXinrandomgeneral} and they proved that 
a shear flow accelerates the propagation \cite{NolenXinrandomshear}, while in dimension $1$ the advection term slows it down \cite{NolenXinrandom1D}. 
Note that, in dimension $1$, it is always possible to turn a reaction-diffusion equation with advection term into a self-adjoint one like (\ref{Cauchy})
by a simple change of variables involving the primitive of the advection term. As we were not able to recover nor to enhance Nolen-Xin's result through this 
change of variables in the present paper, we will not consider any dependence with respect to some advection term. 

Second, propagation phenomenas in random stationary ergodic environment have been investigated by Shigesada and Kawazaki \cite{Shigesadabook}. In their book, they 
considered the case where $a\equiv 1$ and 
$$f(x,\o,s) =\left\{\begin{array}{rcl}  f^+ (s) &\hbox{ if }& X_{2m}(\o)\leq x < X_{2m+1}(\o),\\
                     f^- (s) &\hbox{ if }& X_{2m+1}(\o)\leq x < X_{2m+2}(\o),\\
                    \end{array}\right. \quad \hbox{for some } m\in\mathbb{Z}, $$
with $f^+\geq f^-$ and $(X_n)_n$ is a family of random variables such that $(X_{2m+1}-X_{2m})_m$ and $(X_{2m+2}-X_{2m+1})_m$ are two given families of 
independent and identically distributed variables. Considering the expectations 
$$\ell_1:= \E [X_{2m+1}-X_{2m}] \quad \hbox{and} \quad \ell_2:= \E [X_{2m+2}-X_{2m+1}],$$
they carried some numerical simulations showing that the associated spreading speed is approximately the same as the one associated with 
the $(\ell_1+\ell_2)-$periodic nonlinearity $f_{per} (x,s):= f_+ (s)$ if $x\in [0,\ell_1)$, $f(x,s)=f_-(s)$ if $x \in [\ell_2,\ell_1+\ell_2)$. 
When the variances of $(X_{2m+1}-X_{2m})_m$ and $(X_{2m+2}-X_{2m+1})_m$ increase, only a very slight increase of the spreading speed is observed, which might be due to
numerical remains (see \cite{Shigesadabook}).


\subsection{Dependence results in the periodic framework}

If the dependence between the spreading speed and the coefficients have not been much investigated in random statonary ergodic media, 
it has been extensively studied in the last decade when $a$ and $f$ are periodic with respect to $x$. In this case, another characterization of the speed $w^*$ holds, involving a family of periodic principal eigenvalues associated with the linearization of 
the equation near $u=0$ (see Theorem \ref{thm:BN1} below and the comments afterwards). This characterization thus reduces the dependence of $w^*=w^*(a,f)$ with 
respect to $a$ and $f$ to an eigenvalue optimization problem. 

Note that the case of periodic coefficients is a particular class of random stationary ergodic problems: one just needs to take $\O=\R/L\Z$, where $L$ is the periodic of 
the coefficients, $\mathbb{P}$ the Lesbegue measure on the torus $\R/L\Z$, and $\pi_x \o = \o+y$ (mod L) for all $\o \in \O=\R/L\Z$ and $x\in\R$ in order to turn 
a periodic problem into a random stationary ergodic one. Similarly, it is possible to turn an almost periodic setting into a random strationary ergodic one, but
through a more involved construction \cite{PapanicolaouVaradhan}. 

Many dependence results have been proved using this characterization in the periodic framework. In 2005, Berestycki, Hamel and Nadirashvili \cite{BHNadir} proved that 
if $g\geq f$, then $w^*(a,g)\geq w^*(a,f)$, and if $f$ does not depend on $x$, then $\kappa\mapsto w^*(\kappa a,f)$ is nondecreasing. The second result does not necessarily
hold if $f$ depends on $x$ \cite{dependence}. 
If $f$ takes negative values, then the spreading result still holds under some milder assumptions on the stability of the steady state $0$ \cite{BHNa}. In this case, 
$B\mapsto w^*(1,Bf)$ is increasing if the periodic function $x\mapsto f_s(x,0)$ has a positive average \cite{BHR2}. 

It is possible to extend Theorem \ref{thm:Freidlin} to multidimensional periodic equations \cite{Freidlin} and to space-time periodic coefficients \cite{BHNa}. 
The previous results still hold in the space-time periodic multidimensional media. However, it is not possible to obtain a monotonicity of the 
spreading speed in the diffusion matrix  with respect to the positive matrix ordering \cite{ElSmaily}. Moreover, 
one can compare $w^*(a,f)$ with the spreading associated with 
averaged coefficients in time or in space \cite{dependence}: generally speaking, a high heterogeneity of the coefficients gives a large
spreading speed. 
It is also relevant in multidimensional media to investigate the effect of an incompressible drift term in the equation. Such a drift generally increases 
the spreading speed
\cite{BHNadir}, but the amplitude of this increase depends on the properties of the flow associated with the drift 
(see for instance \cite{Audoly, Largedrift, ElSmailyKirsch, Kiselevadvection, RyzhikZlatos, Zlatosc/A}). 

The main difficulty that one has to face when investigating this type of problems is that the periodic principal eigenvalues associated with the linearization 
near $u=0$ are not necessarily related to a self-adjoint operator. Hence, it could not be expressed as the extremum of some Rayleigh quotient, which makes the
optimization problems for these eigenvalues quite uneasy. This difficulty was overcame by the author in \cite{Steiner}, where a new (non-quadratic) 
integral characterization of the periodic principal eigenvalue of a non-self-adjoint operator has been proved, enabling for example to prove that in dimension $1$, 
if $a\equiv 1$ and $f=f(x,u)$ is periodic in $x$, then taking the Schwarz rearrangement of the growth rate $x\mapsto f_s(x,0)$ increases the associated spreading speed. 
Also, if one considers the rescaling of the reaction rate $f_L (x,u):= f(x/L,u)$ and the diffusion term $a_L (x):= a(x/L)$, then 
$L\mapsto w^*(a_L,f_L)$ is nondecreasing \cite{Steiner}. The limits of this function as $L\to 0^+$ and $L\to +\infty$ were computed respectively in \cite{Steiner}
and \cite{HamelNadinRoques}. This problem was originally stated in a biological modelling framework by Shigesada, Kawazaki and Teramoto  \cite{Shigesada1}, and 
the heuristic interpretation of this monotonicity is that the more the media is fragmented (in the sense that 
the patches of favourable and unfavourable media are both small), the slower the speed of invasion of a species. 
This new integral characterization of periodic principal eigenvalues for non-symmetric operators has been used successfully in several recent papers 
in order to get estimates on the spreading speed \cite{DGM, LLM}.


\subsection{A characterization of the spreading speed involving generalized principal eigenvalues}

The aim of the present paper is to extend such dependence results to random stationary ergodic media. This work was motivated by a recent paper of Berestycki and the author
\cite{BN1}, where a new characterization of the spreading speed in random stationary ergodic media was proved. This formula is very similar to the one which holds in
the periodic framework, except that it involves generalized principal eigenvalues instead of periodic principal ones. 
Namely, define for all $p>0$, $\o\in\O$ and $\phi\in\mathcal{C}^2 (\R)$ the second order elliptic operator 
\begin{equation}
 L_p^\o \phi := \big( a(x,\o)\phi'\big)'-2p a(x,\o) \phi' + \big(p^2 a(x,\o)-pa(x,\o)' + c(x,\o)\big)\phi,
\end{equation}
where we remind to the reader that $c(x,\o):= f_s(x,0)$. As in \cite{BN1}, we associate with this operator two generalized principal eigenvalues:  
\begin{equation}\label{def-kp} \begin{array}{rcrcl}
\overline{k}^\o_p(a,c)&:=& \inf \{&\lambda\in\R, \ \exists  \phi \in\mathcal{A}, \  L_p^\o\phi\leq \lambda\phi&\},\\ 
\underline{k}^\o_p(a,c)&:=& \sup \{&\lambda\in\R, \ \exists  \phi \in\mathcal{A}, \  L_p^\o\phi\geq \lambda\phi&\}, \\ 
\end{array} \end{equation}
where $\mathcal{A}$ is the set of admissible test-functions:
\begin{equation}\label{def-A}
\mathcal{A}:= \big\{\phi \in \mathcal{C}^0(\R), \quad \phi>0 \hbox{ in } \R, \quad \phi'/\phi \in L^\infty (\R), \quad \lim_{|x|\to +\infty} \ln \phi (x)/|x|=0\big\}
\end{equation}
where the inequality should be understood in the viscosity sense. 
Then Berestycki and the author have proved in \cite{BN1} that these two quantities are almost surely deterministic, identical and give a new characterization 
of the spreading speed. 

\begin{thm}\label{thm:BN1}\cite{BN1}

\begin{enumerate}
\item One has $k_p^\o (a,c) \geq k_0^\o (a,c)$ for all $p\in\R$ and $\o\in\O$. 

\item There exists a measurable set $\O_1 \subset \O$, with $\P (\O_1)=1$, such that for all $\o\in\O_1$, 
$\underline{k}_p^\o (a,c)=\overline{k}_p^\o (a,c)$. Moreover, if $\underline{k}_p^\o (a,c) > \Lambda_1^\o (a,c)$, then there exists 
a $\phi\in\mathcal{A}$ such that $L_p^\o \phi = \underline{k}_p^\o(a,c)\phi$ in $\R$. In this case we denote $k_p (a,c):= \underline{k}_p^\o (a,c)=\overline{k}_p^\o (a,c)$.

\item One has $w^*(a,c)=\min_{p>0} \frac{k_p (a,c)}{p}$. 
\end{enumerate}
\end{thm}

This formula is a generalization of the characterization of the spreading speed in periodic media. For general heterogeneous bounded conefficients, this result still holds
except that the two generaliazed principal eigenvalues are not necessarily equal and thus it only gives bounds on the location of the level sets of $u$. 

We are now in position to state our main dependence results. 


\section{Statement of the results}

\subsection{Comparison with homogenized coefficients}

Let first compare the spreading speeds associated with heterogeneous or homogeneous coefficients. It is natural to compare 
the effect of the growth rate $f$ with that of the homogeneous growth rate $\E [f]$. Here and in the sequel (when there is no ambiguity),  
$\mathbb{E} [f] (s):= \mathbb{E} [\o\mapsto f(x,\o,s)]$ and the random stationarity yields that this quantity does not depend on $x$. 
As in periodic media \cite{Steiner}, it is not possible to obtain a comparison with the homogeneous diffusion rate $\E [a]$ and one needs to consider the harmonic average $\E [1/a]^{-1}$ 

\begin{prop} \label{prop-hom}
One has $k_p (a,c)\geq \mathbb{E} [f_s (\cdot,\cdot,0)] + p^2\mathbb{E}[1/a]^{-1}$ for all $p\in\R$ and thus 
$$w^*(a,f) \geq w^* \big(\mathbb{E}[1/a]^{-1},\mathbb{E} [f]\big) = 2\sqrt{\mathbb{E}[1/a]^{-1}\mathbb{E} [f_s(\cdot,\cdot,0)]}.$$
Moreover, if $a$ is constant while $c$ is not, then these inequalities are strict.  
\end{prop}

The second identity is obvious: if the diffusion and the reaction terms are constant in $(x,\o)$, we recover the classical Fisher-KPP equation \cite{KPP}, which spreading speed 
is explicitly determined by the derivative of the reaction term at $s=0$ and the diffusion term. 

In the case of periodic coefficients, this result was proved by Berestycki, Hamel and Roques when $a$ does not depend on $x$ \cite{BHR2}
and by the author in the general framework \cite{Steiner}. 

This result means that heterogeneous coefficients always speed up the propagation. It now remains to investigate what is the amplitude of this acceleration 
with respect to the shape of this heterogeneity. 


\subsection{Dependence with respect to the diffusion}

Consider now the dependence relation between the amplitude of the diffusion term and the spreading speed. 

\begin{prop}\label{prop:diffusion}
 Assume that $c(x,\o)=f_s(x,\o,0)$ does not depend on $(x,\o)$. Then $k_p (\kappa a,c)=\kappa k_p (a,0)+c$ for all $p\in\R$ and thus 
$$\kappa \mapsto w^*(\kappa a,f) \hbox{ is increasing}.$$
\end{prop}

This result is an extension of Theorem 1.10.2. in \cite{BHNadir} and will be proved through similar arguments and Theorem \ref{thm:BN1}. 
It means that an increase of the diffusion rate speeds-up the propagation.
The counter-example exhibited in a periodic framework \cite{dependence} yields that this result cannot be true in general if $c$ depends on $(x,\o)$. 


\subsection{Dependence with respect to the reaction}

\begin{prop}\label{prop:reaction}
1. Assume that $f_s(x,\o,0)\leq g_s(x,\o,0)$ for all $x\in\R$ almost surely in $\o\in\O$. Then 
$k_p \big(a,f_s(\cdot,\cdot,0)\big)\leq k_p \big(a,g_s(\cdot,\cdot,0)\big)$ for all $p\in\R$ and thus  
$$w^*(a,f) \leq w^* (a,g).$$

\smallskip

2. Assume that $f_s(x,\o,0)=f'(0)>0$ does not depend on $(x,\o)$ and that $g=g(x,\o,s)$ is a random stationary ergodic variable satisfying 
(\ref{hyp:KPP}), $g(x,\o,0)=g(x,\o,1)=0$ for all $(x,\o)\in\R\times \O$, and such that $g$ and $g_s(\cdot,\o,0)$ 
are almost surely uniformly continuous and bounded with respect to $x$ uniformly in $s$ and $\mathcal{C}^{1+\gamma}$ with respect 
to $s$ uniformly in $x$. Let 
$$B_*=\sup \{ B\geq 0, f(s)+B g(x,\o,s)>0 \hbox{ for all } (x,s)\in \R \times (0,1) \hbox{ a. e.}\}.$$
Then if $\E [g_s(x,\cdot,0)]\geq 0$, the function $B\mapsto k_p \big(1,B g_s(\cdot,\cdot,0)\big)$ is nondecreasing for all $p\in\R$ and thus, 
if in addition $B_*>0$, 
the function $B\in [0,B_*) \mapsto w^*(1,f+Bg)$ is nondecreasing. 

Moreover, if $c$ is not constant, then these functions are increasing. 
\end{prop}

The heuristical interpretation of 1. is trivial: the higher the growth rate, the faster the propagation. This result was proved in the periodic framework in 
Theorem 1.6 of \cite{BHNadir}.

In 2., we need to restrict to $B\in [0,B_*]$ in order to get a reaction rate $f+Bg$ satisfying the above hypotheses guaranteeing the existence of a spreading speed. 
In periodic media, this result was proved by Berestycki, Hamel and Roques \cite{BHR2}, who gave the following interpretation when the average of 
$g_s(x,0)$ is zero: the larger the oscillations near the average growth rate, the faster the propagation. 
Our result is an extension to random stationary ergodic media. 


\subsection{Dependence with respect to the scaling of the coefficients}

We consider in this section the rescaled coefficients for all $L>0$:
$$a_L (x,\o):= a (x/L, \o) \quad \hbox{and} \quad f_L (x,\o,s) := f(x/L,\o,s) \quad \hbox{for all } (x,\o,s) \in\R\times \O\times [0,1].$$
These coefficients are random stationary ergodic (with measure preserving transformations $\tilde{\pi}_x=\pi_{x/L}$) and satisfy all the hypotheses of Section \ref{sec:hyp}. 
Hence, the spreading speed $w^*(a_L,f_L)$ is well-defined. 

\begin{thm} \label{prop:scaling}
The function $L\mapsto k_p \big(a_L, (f_L)_s(\cdot,\cdot,0)\big)$ is nondecreasing for all $p\in\R$ and thus  $L>0 \mapsto w^* (a_L,f_L)$ is nondecreasing. 
Moreover, if $a$ is constant while $c$ is not, these functions are increasing. 
\end{thm}

This result has been proved in the periodic framework by the author \cite{Steiner}. It means that fragmentation of the environment slows down the propagation 
(see \cite{Steiner, Shigesadabook}). The strict monotonicity when $a$ is constant is a new result, and it is not clear whether it holds when $a$ is 
not constant or not. 


\section{A new formula for the generalized principal eigenvalue of a nonsymmetric operator with random stationary ergodic coefficients}

The aim of this section is to prove the following result:

\begin{thm} \label{formulathm}
One has
\begin{equation} \label{eq:formula} k_p (a,c)=\inf_{\theta\in \mathcal{B}}k_0 (a,c+a |p+\theta|^2),\end{equation}
where
\begin{equation}\label{def-B}\begin{array}{rl}
\mathcal{B}:= \big\{ & \theta : \R\times \O\to \R \hbox{ measurable in } \o\in\O \hbox{ and a.s. bounded in } x\in\R, \\
& \theta (x+y,\o)=\theta (x,\pi_y \o) \hbox{ for all } (x,y,\o)\in\R\times \R\times \O, \ \E [\theta]=0\big\}.\\ \end{array}
\end{equation}
Moreover, for all $p$ in the closure of the set $\{ \tilde{p} \in\R, \ k_{\tilde{p}} (a,c) > k_0 (a,c)\}$,  the infimum in (\ref{eq:formula}) is indeed a minimum. 
\end{thm}

Note that, as $\theta \in\mathcal{B}$, $(x,\o)\mapsto c(x,\o)+a(x,\o) |p+\theta (x,\o)|^2$ is a random stationary ergodic variable and thus 
$k_0 (a,c+a |p+\theta|^2)$ 
is well-defined.

Moreover, it was proved in \cite{BN1} that $\underline{k}_0^\o (a,c) = \overline{k}_0^\o (a,c) = \Lambda_1^\o (a,c)$ almost everywhere, that is
$k_0 (a,c) = \Lambda_1 (a,c)$ since these quantities are almost surely deterministic. Hence one could use $\Lambda_1 (a,c+a |p+\theta|^2)$ instead 
of $k_0 (a,c+a |p+\theta|^2)$ in (\ref{eq:formula}), which will be useful in the sequel. 

This Theorem is a generalization of a similar formula proved by the author in periodic media (see Theorem 2.2 in \cite{Steiner}).
This new characterization of the generalized principal eigenvalue $k_p (a,c)$ will be the key argument in the proofs of 
Propositions \ref{prop-hom}, \ref{prop:diffusion}, \ref{prop:reaction} and Theorem \ref{prop:scaling} below. 


\subsection{Stationarity of the principal eigenfunction}

We first need to improve Theorem \ref{thm:BN1} by showing that, when $k_p^\o (a,c)>\Lambda_1 (a,c)$, then 
the derivative of the logarithm of the associated principal eigenfunction is stationary ergodic with respect to $\o$. We first need to prove the evenness of the eigenvalue with respect to 
$p\in\R$ and the uniqueness up to constant of the eigenfunction. 

\begin{lem}\label{lem:paritek}
For all $p\in\R$, one has 
$$\underline{k}^\o_p (a,c) =\underline{k}^\o_{-p} (a,c) \quad \hbox{almost surely}.$$
\end{lem}

\noindent {\bf Proof.} 
Take $\o\in\O_1$ as in Theorem \ref{thm:BN1} and assume that $\underline{k}_p^\o (a,c) > \Lambda_1^\o (a,c)$ and $\underline{k}_{-p}^\o (a,c)=k_{-p} (a,c) > \Lambda_1^\o (a,c)$. 
Theorem \ref{thm:BN1} gives for all some $\phi\in\mathcal{A}$ and $\psi\in\mathcal{A}$ such that 
$$ \begin{array}{lcll} \big(a(x,\o) \phi'\big)'- 2p a(x,\o) \phi' + \big(c(x,\o)+p^2 a(x,\o)-p a'(x,\o)\big) \phi &=& k_p (a,c) \phi \quad &\hbox{in } \R,\\$$
$$\big( a(x,\o) \psi'\big)'+ 2p a(x,\o) \psi' + \big(c(x,\o)+p^2 a(x,\o)+p a'(x,\o)\big) \psi &=&k_{-p} (a,c) \psi \quad &\hbox{in } \R.\\ \end{array}$$
Multiplying the first equation by $\psi$ and integrating over $(-R,R)$, one gets
$$\begin{array}{cl} &k_p (a,c) \int_{-R}^R \phi\psi dx = \int_{-R}^R \Big(\big( a \phi'\big)'\psi - 2p a \phi'\psi +(c+ap^2 -pa') \phi\psi\Big)dx \\
&\\
= &a(R,\o)  \phi'(R) \psi(R)-a(-R,\o)  \phi'(-R) \psi(-R)-a(R,\o)  \phi(R) \psi'(R)\\
&+a(-R,\o)  \phi(-R) \psi'(-R)- 2p a(R,\o) \phi (R) \psi (R) +2 pa(-R,\o) \phi(-R) \psi (-R)\\
&+\int_{-R}^R \Big( \big(a \psi'\big)'\phi + 2p\phi\big(a \psi\big)' +(c+ap^2 -pa') \phi\psi\Big)dx\\
&\\
=& o\Big( \int_{-R}^R \phi\psi dx\Big) +\int_{-R}^R \Big( \big(a \psi'\big)'\phi + 2pa\phi \psi' +(c+ap^2 +pa') \phi\psi\Big)dx\\
&\\
= &\big( k_{-p} (a,c) +o(1) \big) \int_{-R}^R \phi\psi dx \\
\end{array}$$
using Lemma \ref{lem-cvf} and the equation on $\psi$. Letting $R\to +\infty$ gives the conclusion in this case. 

As $k_p (a,c) \geq k_0 (a,c)$ for all $p\in\R$ and as $p\mapsto k_p (a,c)$ is convex (see \cite{BN1}), one can define 
$$p_+= \sup \{ p\geq 0, \quad k_p (a,c) = k_0 (a,c)\} \quad \hbox{and} \quad p_-= \sup \{ p\geq 0, \quad k_{-p} (a,c) = k_0 (a,c)\}.$$
One has $k_p (a,c) > k_0 (a,c)$ if $p>p_+$ and $k_{-p} (a,c)> k_0 (a,c)$ if $p>p_-$. Hence, the previous step yields 
$k_p (a,c)=k_{-p} (a,c)$ if $p> \max \{p_+,p_-\}$. Assume that $p_+\geq p_-$, then one has $k_{p_+} (a,c) = k_{-p_+} (a,c)=k_0 (a,c)$ by continuity. 
The convexity gives $k_p (a,c)=k_0 (a,c)$ for all $p\in (-p_-,p_+)$. Hence, $p_+=p_-$. This identity can be proved similarly if one assumes $p_-\geq p_+$. 
Eventually, we have proved that $k_p (a,c) = k_{-p} (a,c)=k_0 (a,c)$ if $p\in (-p_+,p_+)$. 
\hfill $\Box$

\begin{lem}\label{lem:unique}
 Let $\O_1$ as in Part 2. of Theorem \ref{thm:BN1}, take $k_p (a,c)>\Lambda_1 (a,c)$, $\o\in\O_1$ and $\phi \in\mathcal{A}$ 
 such that $L_p^\o \phi= k_p (a,c) \phi$ in $\R$. Then such a $\phi$ is unique up to multiplication by a positive constant. 
\end{lem}

When the coefficients are periodic, the uniqueness of the periodic principal eigenvalue is an immediate corollary of the Krein-Rutman theory. 
In the random stationary ergodic setting, this result is new and it is not clear if it holds in the critical framework $k_p (a,c)=\Lambda_1 (a,c)$. 

\bigskip

\noindent {\bf Proof.}
Assume that $\varphi \in\mathcal{A}$ is another solution of $L_p^\o \varphi = k_p (a,c) \varphi$. 
Lemma \ref{lem:paritek} together with Theorem \ref{thm:BN1} yield that there exists $\psi \in\mathcal{A}$ such that $L_{-p}^\o \psi = k_p (a,c) \psi$. 
An easy computation gives $L_p^\o \big( \psi (x) e^{2px}\big) = k_p(a,c) \psi (x) e^{2px}$ for all $x\in\R$. 
Hence, $\phi$ and $x\mapsto \psi (x) e^{2px}$ are two solutions of the same second order ODE, which are not colinear since 
$\lim_{x\to +\infty} \frac{1}{x} \ln\phi (x)=0$ while $\lim_{x\to +\infty} \frac{1}{x} \ln\big(\psi (x)e^{2px}\big)=2p$ and $p\neq 0$ since 
$k_p (a,c)>\Lambda_1 (a,c) = k_0 (a,c)$ by hypothesis. 
Hence, any solution of this ODE is a linear combination of these two functions. In particular, there exists $A,B\in\R$ such that 
$\varphi (x) = A \phi (x) + B \psi (x) e^{2px}$ for all $x\in\R$. But as $\lim_{x\to +\infty} \frac{1}{x} \ln\varphi (x)=0$ since $\varphi \in\mathcal{A}$,
one has $B=0$. \hfill $\Box$

\begin{cor} \label{cor:rse}
 Let $\O_1$ as in Part 2. of Theorem \ref{thm:BN1} and $k_p (a,c)>\Lambda_1 (a,c)$. For all $\o\in\O_1$, let $\phi (0,\o) \in\mathcal{A}$ 
a solution of $L_p^\o \phi (\cdot,\o)= k_p (a,c) \phi (\cdot,\o)$ in $\R$ normalized by $\phi (0,\o)=1$. 
Then $\phi$ is uniquely defined and 
$$\frac{\phi'(x+y,\o)}{\phi (x+y,\o)} = \frac{\phi'(x,\pi_y \o)}{\phi (x,\pi_y\o)} \quad \hbox{ for all } (x,y,\o)\in\R\times \R\times \O_1.$$
\end{cor}

\noindent {\bf Proof.}
The uniqueness immediately follows from Lemma \ref{lem:unique} and the normalization $\phi (0,\o)=1$. 

Next, as $c(x+y,\o)=c(x,\pi_y \o)$ and $a(x+y,\o)=a(x,\pi_y \o)$, the functions $x\mapsto \phi (x+y,\o)$ and $x\mapsto \phi (x,\pi_y\o)$ 
are both principal eigenvalues in $\mathcal{A}$ of the operator $L_p^{\pi_y \o}$. Hence, Lemma \ref{lem:unique} yields that there exists $A\in\R$ (depending on 
$y$ and $\o$)
such that $\phi (x+y,\o) = A \phi (x,\pi_y\o)$. Taking $x=0$ gives $A=\phi (y,\o)$. 
Taking the logarithm of this identity and derivating with respect to $x$, one gets the conclusion. \hfill $\Box$


\subsection{The upper estimate on $k_p (a,c)$}
We start with a preliminary technical lemma.

\begin{lem} \label{lem-cvf}
 Assume that $f\in \mathcal{C}^0 ((0,\infty))$ is a nonnegative function so that $\frac{1}{x} \ln f(x) \to 0$ as $x\to +\infty$ and $f'/f\in L^\infty ((0,\infty))$. 
Then $$\liminf_{R\to +\infty}\frac{f(R)}{\int_0^R f(x)dx} =0.$$
\end{lem}

\noindent {\bf Proof.} We use te same type of arguments as in the proof of Lemma 6.9 in \cite{BerestyckiHamelRossi}. Assume that 
$$\liminf_{R\to +\infty}\Big(f(R)/\int_0^R f(x)dx\Big) >0.$$ Then there exists $\e>0$ and $R_\e>0$ so that 
for all $R\geq R_\e$, one has $f(R)/\int_0^R f(x)dx \geq 2\e$. As $f'/f$ is bounded, there exists $\delta>0$ so that 
$$\forall R\geq R_\e, \ \int_R^{R+\delta} f(x)dx\geq \e\int_0^R f(x)dx.$$ 
Define for all $n\in \N$, $R_n:= R_\e +n\delta$ and $a_n:= \int_{R_{n-1}}^{R_n} f$. As $\cup_{k=1}^n [R_{k-1},R_k]= [R_\e,R_n]\subset (0,R_n]$, one has 
$a_{n+1} \geq \e \sum_{k=1}^n a_k$, from which one can easily deduce 
\begin{equation} \label{eq:an} a_{n+1} \geq \e a_1 (1+\e)^{n-1}.\end{equation}
On the other hand, one has $a_n \leq \delta\|f'/f\|_\infty f(R_n)$. It follows from $\lim_{x\to +\infty} \frac{1}{x} \ln f(x) = 0$ that $$\frac{1}{n}\ln a_n \to 0
\hbox{ as } n\to +\infty,$$ 
which contradicts (\ref{eq:an}).  \hfill $\Box$

\bigskip

We are now in position to prove our upper bound.

\begin{prop} \label{prop:rightineq}
For all $p\in\R$, one has
$$k_p(a,c)\leq\inf_{\theta\in \mathcal{B}}k_0(a,c+ a|p+\theta|^2).$$
\end{prop}

\noindent {\bf Proof of Proposition \ref{prop:rightineq}.} Take $p\in\R$ and $\o\in\O$ such that $\underline{k}_p^\o (a,c)=k_p (a,c)$. 
Let $\epsilon>0$ and $\phi\in \mathcal{A}$ satisfying
$$L_p^\o\phi \geq (k_p (a,c)-\epsilon)\phi \hbox{ in } \R.$$
Take  $\theta \in \mathcal{B}$ and let $\psi(x):=\phi(x) e^{\int_0^x \theta (y,\o)dy}$. As $\phi\in\mathcal{A}$ and $\theta$ is bounded with respect to $x\in\R$ for almost every $\o\in\O$, 
the function $\psi'/\psi = \phi'/\phi + \theta$ is bounded in $x$, and as $\theta\in\mathcal{B}$ and $\phi\in\mathcal{A}$, one has 
$\lim_{|x|\to +\infty} \frac{1}{|x|}\ln \psi (x)=0$ almost surely since the Birkhoff ergodic theorem yields 
$\frac{1}{|x|}\int_0^x \theta (y,\o)dy \to \E [\theta (x,\cdot)]=0$ as $|x|\to +\infty$ almost surely. 
Hence, $\psi \in\mathcal{A}$ almost surely. 
Moreover, one has 
$$(a\psi')'-2a(p+\theta)\psi' +\big(c+a\theta^2+2pa \theta+p^2 a- (a(\theta+p))'\big)\psi\geq (k_p (a,c)-\epsilon)\psi \hbox{ in } \R.$$
Multiplying by $\psi$ and integrating over $x\in (-R,R)$, Lemma \ref{lem-cvf} gives
$$\begin{array}{cl}
& (k_p (a,c)-\epsilon)\int_{-R}^R\psi^2 dx\\
&\leq \int_{-R}^R\Big( \big(a\psi'\big)'\psi-2 a(p+\theta \big)\psi'\psi +\big(c+a|\theta+p|^2- (a(\theta+p))'\big)\psi^2\Big)dx \\
&\\
\leq & a(R,\o)\psi'(R) \psi (R) - a(-R,\o)\psi'(-R)\psi(-R) \\
&-a(R,\o)\big(p+\theta (R)\big)\psi^2 (R)+a(-R,\o)\big(p+\theta (-R)\big)\psi^2 (-R) \\
&+ \int_{-R}^R\Big( -a(\psi')^2+\big(c+a|\theta+p|^2\big)\psi^2\Big)dx\\
&\\
\leq & 2\big(\|a\|_\infty \|\psi'/\psi\|_\infty +(p+\|\theta\|)\|a \|\big) \big(\psi^2 (R)+\psi^2 (-R)\big) 
+\Lambda_1^\omega \big(a,c+a|p+\theta|^2\big)\int_{-R}^R \psi^2 dx\\
&\\
\leq & \Big( o(1)+ \Lambda_1^\omega \big(a,c+a|p+\theta|^2\big)\Big)\int_{-R}^R \psi^2 dx\\
\end{array}$$
where $\Lambda_1^\omega$ was defined in (\ref{def:Lambda1}). 
Hence: 
$k_p (a,c)-\epsilon \leq\Lambda_1^\omega \big(a ,c+a|p+\theta|^2\big)$ almost surely,
for all $\epsilon>0$ and thus 
$$k_p (a,c) \leq\Lambda_1^\omega \big(a ,c+a|p+\theta|^2\big) \quad \hbox{ almost surely}.$$

Consider now a sequence $(\theta_n)_n$ in $\mathcal{B}$ such that 
$$\lim_{n\to +\infty} k_0 (a,c+a|p+\theta_n|^2) = \inf_{\theta \in\mathcal{B}} k_0 (a,c+a|p+\theta|^2).$$
for all $n$, there exists $\O^n\subset \P$ of measure $1$ such that 
$$\underline{k}_0^\o (a,c+a|p+\theta_n|^2) = \overline{k}_0^\o(a,c+a|p+\theta_n|^2) = k_0 (a,c+a|p+\theta_n|^2) \quad \hbox{and}$$
$$\Lambda_1^\o(a,c+a|p+\theta_n|^2) = \Lambda_1(a,c+a|p+\theta_n|^2)$$ 
for all $\o \in \O^n$. Let $\O_1:= \cap_{n} \O^n$, then 
as $k_0 \equiv \Lambda_1$ (see Lemma 5.1 of \cite{BN1}), 
for all $\o\in\O_1$, one has
$$k_p (a,c) \leq\Lambda_1^\omega \big(a ,c+a|p+\theta_n|^2\big) = k_0  \big(a ,c+a|p+\theta_n|^2\big) \to  \inf_{\theta \in\mathcal{B}} k_0 (a,c+a|p+\theta|^2) \quad \hbox{ as } n\to +\infty.$$

\hfill $\Box$


\subsection{The lower estimate on $k_p (a,c)$}

\begin{prop} \label{prop:leftineq}
One has
\begin{equation}\inf_{\theta\in \mathcal{B}}k_0 (a,c+a|p+\theta|^2)\leq k_p (a,c).\end{equation}
\end{prop}

\noindent {\bf Proof of Proposition \ref{prop:leftineq}.} Take $p\in\R$, $\o\in\O_1$ as in Theorem \ref{thm:BN1} and such that $k_0 (a,c) = \Lambda_1^\o (a,c)$.
Assume that $k_p (a,c)> k_0 (a,c)$.
Take  and $\phi\in\mathcal{A}$ such that 
$$ \big(a(x,\o)\phi'\big)'-2p a(x,\o)\phi' +\big(c(x,\o)+p^2 a(x,\o)-pa'(x,\o)\big)\phi=k_p (a,c)\phi \hbox{ in } \R.$$
Similarly, Lemma \ref{lem:paritek} yields that one can assume the existence of $\psi\in\mathcal{A}$ satisfying
$$ \big(a(x,\o)\psi'\big)'+2p a(x,\o)\psi' +\big(c(x,\o)+p^2 a(x,\o)+pa'(x,\o)\big)\psi=k_p (a,c)\psi \hbox{ in } \R.$$
Let $\alpha := \sqrt{\phi\psi}\in\mathcal{A}$ and 
$$\theta:=-\displaystyle\frac{\phi'}{2\phi} +\displaystyle\frac{\psi'}{2\psi}.$$
Corollary \ref{cor:rse} yields that 
$$\theta (x+y,\o) = \theta (x,\pi_y \o) \quad \hbox{ for all } (x,y,\o)\in\R\times \R\times \O_1.$$
Hence $\theta$ is a random stationary ergodic funtion and thus the Birkhoff ergodic theorem yields that for all $\o\in\O_1$, 
$$\E [\theta] = \lim_{x\to +\infty} \int_0^x \theta (y,\o)dy.$$
But as $\phi,\psi \in\mathcal{A}$, the right hand side is $0$. Hence $\E[\theta]=0$
 and thus $\theta \in\mathcal{B}$. 

Next, we compute $\alpha' = \displaystyle\frac{\phi'\psi + \psi'\phi}{2\sqrt{\phi\psi}}$ and:
$$ \begin{array}{rcl}
\big( a\alpha'\big)' &=& \displaystyle\frac{\big(a\phi'\big)'\psi + \big( a\psi'\big)'\phi+2 a\phi'\psi'}{2\sqrt{\phi\psi}}-\frac{a(\phi'\psi+\psi'\phi)^2}{4 (\phi\psi)^{3/2}}\\ 
&&\\
&=& \displaystyle\frac{\big(a\phi'\big)'\psi + \big( a\psi'\big)'\phi}{2\sqrt{\phi\psi}}-\frac{a(\phi'\psi-\psi'\phi)^2}{4 (\phi\psi)^{3/2}}\\ 
&&\\
&=& \displaystyle\frac{\big(a\phi'\big)'\psi + \big( a\psi'\big)'\phi}{2\sqrt{\phi\psi}}-a \theta^2\alpha\\ 
&&\\
&=& \displaystyle\frac{2p a\phi'\psi+pa'\phi\psi - 2pa\psi'\phi-pa'\psi\phi}{2\sqrt{\phi\psi}}-\big( c+p^2a -k_p (a,c) +\theta^2a\big)\alpha\\ 
&&\\
&= & \displaystyle\frac{2pa\big(\phi'\psi -  \psi'\phi\big)}{2\sqrt{\phi\psi}}-\big( c+p^2a -k_p (a,c) + \theta^2a\big)\alpha\\ 
&&\\
&=& -\big( c+|p+\theta|^2 a-k_p (a,c) \big)\alpha.\\ 
\end{array}$$

In other words, we have constructed an exact principal eigenfunction $\alpha \in\mathcal{A}$. Proposition 2.2 in \cite{BN1} thus gives
$$\underline{k}_0^\o\big(a,c+|p+\theta|^2 a\big)=\overline{k}_0^\o\big(a,c+|p+\theta|^2 a\big)=k_p (a,c).$$
On the other hand, Proposition \ref{prop:rightineq} yields 
$$k_p (a,c)\leq \inf_{\tilde{\theta} \in\mathcal{B}} k_0\big(a,c+|p+\tilde{\theta}|^2 a\big).$$
Hence $\theta$ minimizes this infimum. 

Next, as in the proof of Lemma \ref{lem:paritek}, we know that there exists $p_+\geq 0$ such that 
$k_p (a,c)=k_0 (a,c)$ if $p\in [0,p_+]$ while $k_p (a,c)>k_0 (a,c)$ if $p>p_+$. 
Take $p>p_+$, let $t=p_+/p \in [0,1)$ and $\theta\in\mathcal{B}$ such that $k_0\big(a,c+|p+\theta|^2 a\big) = k_p (a,c)$. 
As $c\mapsto k_0 (a,c)$ is convex, one has 
$$\begin{array}{rcl} k_0 (a,c+|p_+ +t\theta|^2 a)&=& k_0(a,c+t^2|p +\theta|^2 a)\\
&&\\
   &\leq &(1-t^2) k_0 (a,c)+t^2 k_0(a,c+|p+\theta|^2 a) \\
&&\\
&\leq& (1-t^2) k_0(a,c)+t^2 k_p (a,c).\\
  \end{array}$$

Hence, $\inf_{\theta \in\mathcal{B}} k_0 (a,c+|p_+ +\theta|^2 a) \leq (1-t^2) k_0(a,c)+t^2 k_p (a,c)$. 
Take a sequence $(p_n)_n$ such that $p_n>p_+$ for all $n$ and $\lim_{n\to +\infty} p_n = p_+$, 
we have thus proved that $\inf_{\theta \in\mathcal{B}} k_0 (a,c+|p_+ +\theta|^2 a) \leq (1-t_n^2) k_0 (a,c)+t_n^2 k_{p_n} (a,c)$
for all $\o\in \cap_{n}\O_n$, with $t_n = p_+/p_n$. Letting $n\to +\infty$, 
as $k_{p_n}(a,c) \to k_{p_+} (a,c)$ as $n\to +\infty$ by convexity, we eventually get
$$\inf_{\theta \in\mathcal{B}} k_0 (a,c+|p_+ +\theta|^2 a) \leq k_{p_+} (a,c)=k_0 (a,c).$$
If $p\in (0,p_+)$, then letting $t= p/p_+$, one gets 
$$\inf_{\theta \in\mathcal{B}} k_0 (a,c+|p +\theta|^2 a) \leq k_0(a,c+|p_+ +t\theta|^2 a) 
\leq (1-t^2) k_0 (a,c)+t^2 k_{p_+} (a,c)= k_0 (a,c)=k_p (a,c)$$
almost surely, which concludes the proof. 
\hfill $\Box$

\begin{lem}\label{lem:minimizer}
 For all $p$ in the closure of the set $\{ \tilde{p} \in\R, \ k_{\tilde{p}} (a,c) > k_0 (a,c)\}$,  the infimum in (\ref{eq:formula}) is indeed a minimum. 
\end{lem}

Note that we have already noticed that, by convexity, $$\{ \tilde{p} \in\R, \ k_{\tilde{p}} (a,c) > k_0 (a,c)\} = (-\infty,-p_-)\cup (p_+,+\infty)$$
but 
such an accurate description will not be needed in the proof of this Lemma. 

\bigskip

\noindent {\bf Proof.}
If $k_p (a,c) > k_0 (a,c)$, then this result was part of the proof of Proposition \ref{prop:leftineq}. 

Consider now a sequence $(p_n)_n$ converging to $p_\infty \in \R$ such that $k_{p_n} (a,c) > k_0 (a,c)$ for all $n\in\N$ and 
$\big(k_{p_n}(a,c)\big)_n$ is a decreasing sequence. Let $\theta_n \in\mathcal{B}$ 
such that $k_{p_n} (a,c) = k_0 (a,c+|p_n+\theta_n|^2)$ for all $n$. 
Then Lemma \ref{lem-compint}, that will be proved later without making use of Lemma \ref{lem:minimizer}, yields 
$$\E [ \tilde{c}+\tilde{a} |p_n +\tilde{\theta_n}|^2] \leq k_{p_n} (a,c),$$
where $\tilde{c} (\o):= c(0,\o)$, $\tilde{a} (\o):= a(0,\o)$ and $\tilde{\theta_n} (\o):= \theta_n(0,\o)$. Hence, $\tilde{\theta_n}$ is bounded in $L^2(\O)$
and we could thus assume that it converges weakly in $L^2 (\O)$ to some $\tilde{\theta_\infty}$. 
Mazur's lemma yields that a sequence of convex combinations $(\xi_k)_k$ of $(\tilde{\theta_n})_{k\geq n}$ converges to $\tilde{\theta_\infty}$ strongly 
in $L^2(\O)$. 
It is known (see \cite{Oleinik}) that there exists a set of full measure $\O_2$ such that for all $\o\in \O_2$, 
the sequence of functions $x\mapsto \xi_n (x,\o) = \xi_n (\pi_x \o)$ converges to $\theta_\infty (x,\o):= \tilde{\theta_\infty} (\pi_x \o)$ in 
$L^2_{loc} (\R)$. We could assume that $\O_1 \subset \O_2$ by changing $\O_1$ if necessary. 

Coming back to the definition of $\theta_n$ in the proof of Proposition \ref{prop:leftineq}, we get 
$\theta_n = -\frac{\phi_n'}{2\phi_n}+\frac{\psi_n'}{2\psi_n}$, where $\phi_n$ and $\psi_n$ are eigenvalues associated with $L_{p_n}^\o$ and $L_{-p_n}^\o$. 
As $a$, $a'$, $1/a$ and $c$ are uniformly bounded and continuous almost surely and the sequence $(p_n)_n$ is bounded, it follows from the elliptic Harnack inequality that 
$\phi_n'/\phi_n$ and $\psi_n'/\psi_n$ are bounded and that this bound does not depend on $n$. Hence, $(\theta_n (\cdot,\o))_n$ is uniformly bounded
over $\R$ for all $\o\in \O_1$. Thus, $\theta_\infty$ is almost surely bounded. 
As it is clear that $\E [\theta_\infty] = 0$ by weak convergence, we conclude that $\theta_\infty \in\mathcal{B}$. 

Next, we know that 
$$\begin{array}{rcl} k_{p_n} (a,c)&=&k_0 (a,c+a|p_n+\theta_n|^2) = \Lambda_1^\o (a,c+a|p_n+\theta_n|^2) \\
   &&\\
   &=& \sup_{\alpha \in H^1 (\R)\backslash \{0\}} \displaystyle \frac{\int_\R \Big(-a(x,\o)(\alpha'(x))^2 + \big( c(x,\o)+ a(x,\o) |p_n +\theta_n (x,\o)|^2 \big)\alpha^2 (x)\Big)dx}{\int_\R \alpha^2 (x)dx}.\\
  \end{array}$$

  By convexity and monotonicity of the sequence $\big(k_{p_n}(a,c)\big)_n$, we get 
  $$\sup_{\alpha \in H^1 (\R)\backslash \{0\}} \displaystyle \frac{\int_\R \Big(-a(x,\o)(\alpha'(x))^2 + \big( c(x,\o)+ a(x,\o) |p_n +\xi_n (x,\o)|^2 \big)\alpha^2 (x)\Big)dx}{\int_\R \alpha^2 (x)dx}
  \leq k_{p_n}(a,c).$$
  
Take any smooth compactly supported function $\alpha \not\equiv 0$. 
Then for all $\o\in\O_1$,  
$$\lim_{n\to +\infty} \int_\R a(x,\o)|p_n +\xi_n (x,\o)|^2 \alpha (x)^2 dx \geq \int_\R a(x,\o)|p_\infty +\theta_\infty (x,\o)|^2 \alpha (x)^2 dx$$
by weak convergence. Hence, as smooth compactly supported functions are dense in $H^1 (\R)$ and $p\mapsto k_p (a,c)$ is convex and thus continuous, we get 
$$\begin{array}{rcl} k_{p_\infty} (a,c)&=& \lim_{n\to +\infty} k_{p_n} (a,c)\\
&&\\
&\geq& \sup_{\alpha \in H^1 (\R)\backslash \{0\}} 
\displaystyle \frac{\int_\R \Big(-a(x,\o)(\alpha'(x))^2 + \big( c(x,\o)+ a(x,\o) |p_\infty +\theta_\infty (x,\o)|^2 \big)\alpha^2 (x)\Big)dx}{\int_\R \alpha^2 (x)dx}\\
&&\\
&=& \Lambda_1^\o (a,c+a|p_\infty +\theta_\infty|^2).\\ 
\end{array}$$
Hence, by Proposition \ref{prop:leftineq}, $\theta_\infty$ is an admissible minimizer, which concludes the proof. 
\hfill $\Box$

\bigskip

\noindent {\bf Proof of Theorem \ref{formulathm}} It immediately follows from Propositions \ref{prop:rightineq}, \ref{prop:leftineq} and Lemma \ref{lem:minimizer}. \hfill $\Box$


\section{Proof of the dependence results}


\subsection{Proof of the comparison with respect to homogenized coefficients}

\begin{lem} \label{lem-compint}
 One has $k_0 (a,c) \geq \E [c]$. 
\end{lem}

\noindent {\bf Proof of Lemma \ref{lem-compint}.}
Take $\O_1$ as in Theorem \ref{thm:BN1} and $\lambda > \overline{k}_0^\o (a,c)$. The definition of $\overline{k}_0^\o$ yields that there exists $\phi\in\mathcal{A}$ such that 
$$\big( a(x,\o)\phi'\big)'  +c(x,\o)\phi \leq \lambda \phi \quad \hbox{in } \R.$$
Dividing this inequality by $\phi$ and integrating by parts over $(-R,R)$, one gets 
$$\begin{array}{rcl}2R\lambda &\geq&  \displaystyle\int_{-R}^R \frac{\big( a(x,\o)\phi'\big)'(x)}{\phi (x)} dx+\int_{-R}^R c(x,\o)dx\\
   &\geq& a(R,\o) \displaystyle\frac{\phi'(R)}{\phi (R)} - a(-R,\o) \displaystyle\frac{\phi'(-R)}{\phi (-R)}
+\int_{-R}^R \displaystyle\frac{a(x,\o)\big(\phi'(x)\big)^2}{\phi (x)^2} dx 
+ \int_{-R}^R  c(x,\o)dx\\
&\geq & a(R,\o) \displaystyle\frac{\phi'(R)}{\phi (R)} - a(-R,\o) \displaystyle\frac{\phi'(-R)}{\phi (-R)}+ \int_{-R}^R  c(x,\o)dx.\\
  \end{array}$$
The Birkhoff ergodic theorem yields that $\lim_{R\to +\infty} \frac{1}{2R} \int_{-R}^R c(x,\o)dx = \E [c]$ almost surely. 
Dividing the above set of inequalities by $2R$ and letting $R\to +\infty$, as $a$ and $\phi'/\phi$ are bounded, one eventually gets 
$\lambda \geq \E [c]$ almost surely. As this inequality is true for any $\lambda \geq \overline{k}^\o_0 (a,c)$ and as $k_0 (a,c) = \overline{k}^\o_0 (a,c)$ almost surely, this gives the conclusion. \hfill $\Box$

\bigskip

\noindent {\bf Proof of Proposition \ref{prop-hom}.}
For almost every $\o\in\O$, Theorem \ref{formulathm} yields that there exists a sequence $(\theta_n)_n$ in $\mathcal{B}$ such that 
$$\underline{k}_0^\o (a,a|p+\theta_n|^2 +c) \leq \underline{k}_p^\o (a,c)+1/n.$$
On the other hand, we know from Lemma \ref{lem-compint} that for all $n$, there exists $\O_n\subset \O$ such that $\P (\O_n)=1$ and 
$$\underline{k}^\o_0 (a,a|p+\theta_n|^2 +c) \geq \mathbb{E} [\tilde{a}|p+\theta_n (x,\cdot)|^2+\tilde{c}] \quad \hbox{ for all } \o\in\O_n,$$
where we denote $a(x,\o)= \tilde{a} (\tau_x \o)$ and $c(x,\o)= \tilde{c} (\tau_x \o)$ (which is always possible since $a$ and $c$ are a random stationary ergodic variables).
Hence, for all $\o\in\cap_n\O_n$, that is, for almost every $\o\in\O$, and for all $n$, one has 
$$\underline{k}_p^\o (a,c)+1/n\geq \mathbb{E} [\tilde{a}|p+\theta_n (x,\cdot)|^2+\tilde{c}] \geq \inf_{\theta \in\mathcal{B}}\mathbb{E} [\tilde{a}|p+\theta (x,\cdot)|^2] + \mathbb{E}[\tilde{c}],$$
and one can pass to the limit $n\to +\infty$ in order to get rid of the $1/n$ in the left-hand side. 

Let $\tilde{\mathcal{B}} = \{ \tilde{\theta} \in L^2 (\O), \ \mathbb{E} [\tilde{\theta}]=0\}$. It is easily checked that 
for all $x\in\R$, 
$$\{ \theta (x,\cdot), \ \theta \in\mathcal{B}\} = \{ \tilde{\theta} \in L^\infty (\O), \ \mathbb{E} [\tilde{\theta}]=0\}$$ 
(just define $\theta (x,\o)= \tilde{\theta}(\tau_x \o)$) and that this set is dense in $\tilde{B}$. 
Hence, for almost every $\o\in\O$:
$$\underline{k}_p^\o (a,c)\geq \inf_{\tilde{\theta} \in \tilde{\mathcal{B}}}\mathbb{E} [\tilde{a}|p+\tilde{\theta}|^2] + \mathbb{E}[\tilde{c}].$$
Next, the function $\tilde{\theta}\in \tilde{\mathcal{B}}\mapsto\mathbb{E} [\tilde{a}|p+\tilde{\theta}|^2] $ is convex and if $\tilde{\theta}_0$ is a critical point of this function, then $\mathbb{E} [\tilde{a}(\tilde{\theta}_0+p)h]=0$ for all $h\in\tilde{\mathcal{B}}$. But as $\tilde{\mathcal{B}}^\perp$ is exactly the set of constant functions over $\O$, it would then follow that $\tilde{a}(\o)\big(\tilde{\theta}_0(\o)+p\big)=C$ a.e. for some $C\in\R$.
As $\mathbb{E}[\tilde{\theta}_0]=0$, the constant $C$ would be $C=p/\mathbb{E}[1/a]$. It easily follow from this analysis that 
$\tilde{\theta}_0 = p\Big( \displaystyle\frac{1}{\mathbb{E}[1/\tilde{a}]\tilde{a}}-1\Big)$ is the unique minimizer of 
 the function $\tilde{\theta}\in \tilde{\mathcal{B}}\mapsto\mathbb{E} [\tilde{a}|p+\tilde{\theta}|^2]$. 
Thus, 
$$\inf_{\tilde{\theta} \in \tilde{\mathcal{B}}}\mathbb{E} [\tilde{a}|p+\tilde{\theta}|^2]
= \mathbb{E} \Big[ \tilde{a} \times \frac{p^2}{\tilde{a}^2 \mathbb{E}[1/\tilde{a}]^2}\Big] = \frac{p^2}{\mathbb{E}[1/\tilde{a}]}.$$ 
 We conclude that 
$$\underline{k}_p^\o (a,c)\geq \mathbb{E} [\tilde{c}] + \frac{p^2}{\mathbb{E}[1/\tilde{a}]}$$
almost surely and thus 
$$w^*(a,c) = \min_{p>0} \frac{k_p (a,c)}{p} \geq \min_{p>0} \Big(  \frac{\mathbb{E} [\tilde{c}]}{p} + \frac{p}{\mathbb{E}[1/\tilde{a}]}\Big)=2 \sqrt{\frac{\mathbb{E}[\tilde{c}]}{ \mathbb{E}[1/\tilde{a}]}}.$$

Next, assume that $a$ is constant while $c$ is not. Up to some rescaling we could assume that $a\equiv 1$. 
Define $p_+\geq 0$ as in the proof of Proposition \ref{prop:leftineq} and take $p\geq p_+$. 
Lemma \ref{lem:minimizer} thus yields that there exists $\theta \in\mathcal{B}$ such that $k_{p} (1,c) = k_0 (1,c+ |p +\theta|^2)$. 
Assume by contradiction that $k_p (1,c) = \E [\tilde{c}]+p^2$. Then all the previous inequalities are equalities and, 
as $\tilde{\theta}_0$ is the unique minimizer of $\tilde{\theta}\in \tilde{\mathcal{B}}\mapsto\mathbb{E} [|p+\tilde{\theta}|^2]$, it would then 
follow that $\theta (0,\o)=\tilde{\theta}_0 (\o)$ a.e. But here as $\tilde{\theta}_0 = p\Big( \displaystyle\frac{1}{\mathbb{E}[1/\tilde{a}]\tilde{a}}-1\Big)$
and $\tilde{a}$ is constant, one has $\tilde{\theta}_0 \equiv 0$. Hence $\theta\equiv 0$ and the definition of $\theta$ gives
$k_p (a,c) = k_0 (1, c+p^2)=k_0 (1,c)+p^2$.
Putting this into equation $L_p^\o \phi = k_p (1,c)\phi$, we get for all $\o\in\O_1$:
$$\phi'' -2p \phi' + c(x,\o) \phi = k_0 (1,c) \phi \quad \hbox{ in } \R,$$
where $\phi = \phi (x,\o)$ is defined as in Corollary \ref{cor:rse}. Dividing by $\phi$ and integrating, we get 
$$2 k_0 (1,c) R = \displaystyle\frac{\phi'(R)}{\phi (R)} -  \displaystyle\frac{\phi'(-R)}{\phi (-R)}
+\int_{-R}^R \displaystyle\frac{\big(\phi'(x)\big)^2}{\phi (x)^2} dx +\int_{-R}^R c(x,\o)dx -2p \ln\phi (R,\o)+2p \ln \phi (R,\o).$$
Hence, diving by $2R$, letting $R\to +\infty$ and using Birkhoff ergodic theorem, $\phi\in\mathcal{A}$ and the fact that $\phi'/\phi$, and thus $(\phi'/\phi)^2$, 
is a bounded random stationary
ergodic function by Corollary \ref{cor:rse}, we obtain 
$$ k_0 (1,c)  = \E[(\phi'/\phi)^2] +\E[c].$$
But we have assumed by contradiction that $k_p (1,c) = \E [\tilde{c}]+p^2$ and we have obtained $k_p (1,c)=k_0 (1,c)+p^2$. We thus eventually get
$\E[(\phi'/\phi)^2]=0 $, that is, $\displaystyle\frac{\phi'}{\phi} (0,\o)= 0$ a.e. and thus $\phi'/\phi \equiv 0$ by stationarity, meaning that $\phi$ is a constant. 
Lastly, coming back to the equation satisfied by $\phi$, this would imply that $c$ is a constant, a contradiction. Hence
$$k_p (1,c) > \E [\tilde{c}]+p^2.$$
Now, if $p\in [0,p_+)$, one has $k_p (1,c) = k_{p_+}(1,c) > \E [\tilde{c}]+p_+^2\geq \E [\tilde{c}]+p^2$ and thus the strict inequality still holds. 
The strict inequality on $w^*(1,c)$ easily follows. 
\hfill $\Box$


\subsection{Proof of the monotonicity with respect to the diffusion}

\noindent {\bf Proof of Proposition \ref{prop:diffusion}.}
As $c$ does not depend on $x$, one has
\begin{equation} \label{eq:k_p(kappa a,c)} k_p (\kappa a,c) = k_p (\kappa a,0)+c = \kappa k_p (a,0)+c \quad \hbox{a.e.}\end{equation}
On the other hand, we have already noticed that $k_p (a,0) \geq p^2 \E [1/a]^{-1}$, hence $k_p (a,0)>0$, for all $p\in\R$. 
Thus $\kappa \mapsto k_p (\kappa a,c)$ is increasing from (\ref{eq:k_p(kappa a,c)})
and the conclusion follows from part 3. of Theorem \ref{thm:BN1} since $w^*(\kappa a,c)$ is a minimum. \hfill $\Box$


\subsection{Proof of the monotonicity with respect to the reaction}

\noindent {\bf Proof of Proposition \ref{prop:reaction}.}
1. As $c(x,\o)=f_s(x,\o,0) \leq g_s(x,\o,0)=: d(x,\o)$ for all $x\in\R$ almost everywhere, it immediately follows from (\ref{def-kp}) that 
$k_p (a,c)\leq k_p (a,d)$ almost everywhere. The conclusion follows from part 3. of Theorem \ref{thm:BN1}. 

\smallskip

2. Let $c(x,\o):= g_s(x,\o,0)$ and $F(B):= k_p (1,Bc)$. When $c$ is periodic in $x$, the result was proved in Proposition 4.8 of \cite{BHR2} 
by computing the derivative of $F$ at $B=0$. 
We did not manage to compute such a derivative for random stationary coefficients and we will thus use a different argument.

We know from the same types of arguments as in 
Proposition 3.6 in \cite{BHNa} that $B\mapsto \overline{k}_p^\o (1,Bc)$ is a convex function
for all $\o\in\O$. Hence $F$ is convex.  
 Indeed, 
we know from the proof of Proposition \ref{prop-hom} that $F(B) = k_p (1,Bc) \geq B \E[c]+p^2\geq p^2=F(0)$ since $\E [c]\geq 0$ by hypothesis. 
Hence, the monotonicty follows from the convexity and the inequality $F(B)\geq F(0)$ for all $B\geq 0$. 
If $c$ is not constant, then Proposition \ref{prop-hom} yields that $F(B)> F(0)$ for all $B>0$ and the strict monotonicity follows by convexity. 

As $F(B)=k_p (1,Bc)$ is nondecreasing and $k_p \big(1,Bc+f'(0)\big)=k_p (1,Bc)+f'(0)$ since $f'(0)$ is a constant, 
$w^*(1,B g + f)= \min_{p>0} \frac{1}{p} k_p \big(1,Bc+f'(0)\big)$
is also nondecreasing with respect to $B$, and increasing if $c$ is not constant. \hfill $\Box$


\subsection{Proof of the monotonicity with respect to the scaling}

\begin{lem}\label{lem:scaling}
 For all $L>0$, one has $k_{p} (a_L,f_L) = \frac{1}{L^2}k_{pL} (a,L^2f)$ for all $p\in\R$. 
\end{lem}

\noindent {\bf Proof.} Take $p \in\R$ and $\o\in\O$. 
Assume that $\phi \in \mathcal{A}$ and $\lambda\in\R$ satisfy 
$$\big( a_L(x,\o)\phi'\big)' -2p a_L (x,\o) \phi' + \big( p^2 a_L (x,\o)-p a_L'(x,\o) + c_L(x,\o) \big) \phi \geq \lambda \phi \hbox{ in } \R.$$
Let $\psi (x) := \phi (Lx)$. A straightforward computation gives
$$\frac{1}{L^2}\big( a(x,\o)\psi'\big)' -\frac{2p}{L} a (x,\o) \psi' + \big( p^2 a (x,\o) -\frac{p}{L} a'(x,\o) + c(x,\o) \big) \psi \geq \lambda \psi \hbox{ in } \R,$$
which can also be written 
$$\big( a(x,\o)\psi'\big)' -2pL a (x,\o) \psi' + \big( (Lp)^2 a (x,\o) -L p a'(x,\o) + L^2 c(x,\o) \big) \psi \geq L^2\lambda \psi \hbox{ in } \R,$$
It follows from (\ref{def-kp}) that $\frac{1}{L^2}\underline{k}^\o_{Lp} (a,L^2c) \leq \underline{k}^\o_{p} (a_L,c_L)$. Similarly, one can prove that 
$\frac{1}{L^2}\overline{k}^\o_{Lp} (a,L^2c) \geq \overline{k}^\o_{p} (a_L,c_L)$. 
Hence, considering $\o\in\O$ such that $k_p (a_L,c_L) = \overline{k}^\o_{p} (a_L,c_L) = \underline{k}^\o_{p} (a_L,c_L)$ and 
$k_{Lp}(a,L^2c) =\overline{k}^\o_{Lp} (a,L^2c) =\underline{k}^\o_{Lp} (a,L^2c)$, we get the required identity. 
\hfill $\Box$

\bigskip

\noindent {\bf Proof of Theorem \ref{prop:scaling}.} Let $L>1$. 
It follows from Theorem \ref{formulathm} that 
$$k_{Lp} (a,L^2 c)=\inf_{\theta\in \mathcal{B}}\underline{k}_0^\o \big(a,a(\theta +Lp)^2+L^2c \big) \quad \hbox{almost surely}.$$
Letting $\zeta=\theta/L$, one gets
$$\begin{array}{rcl}
k_{Lp} (a,L^2 c)&=& \inf_{\zeta\in \mathcal{B}}\underline{k}_0^\o\Big(a,L^2 \big(a(\zeta+p)^2+c \big)\Big)\\
&&\\
&\geq & L^2 \inf_{\zeta\in \mathcal{B}}\underline{k}_0^\o\big(a,a(\zeta+p)^2+c \big)\\
&&\\
&=& L^2 k_p (a,c) \\ \end{array}$$
almost surely since $L\mapsto k_0 \big(a,L^2 d\big)$ is convex for all bounded uniformly continuous function $d$ (see Proposition 3.6 in \cite{BHNa}). 
It follows from Lemma \ref{lem:scaling} that 
$$ \begin{array}{rcl} w^* (a_L, f_L) &=& \displaystyle\min_{p>0} \frac{k_p (a_L,c_L)}{p}=\displaystyle\min_{p>0} \frac{k_{Lp} (a,L^2 c)}{L^2p}\\
&&\\
 &\geq & \displaystyle\min_{p>0}\frac{k_p(a, c)}{p} \\
&&\\
&=& w^*(a,f).\\ \end{array}$$
As $a$ and $f$ are arbitrary in all this computation, the monotonicity immediately follows. 

If $a$ is a constant while $c$ is not, then the same arguments as in the proof of Proposition \ref{prop:reaction} yield
$$\underline{k}_0^\o\Big(a,L^2 \big(a(\zeta+p)^2+c \big)\Big)> L^2 \underline{k}_0^\o\big(a,a(\zeta+p)^2+c \big)$$
when $L>1$ and the strict monotonicity follows.
\hfill $\Box$


\end{document}